\documentclass[12pt]{article}
\usepackage[margin=3cm]{geometry}
\usepackage{amssymb}

\newtheorem{thm}{Theorem}
\newtheorem{prop}[thm]{Proposition}
\newtheorem{lem}[thm]{Lemma}
\newtheorem{cor}[thm]{Corollary}

\newcommand{\pf}{\noindent {\bf Proof:}\ }
\newcommand{\fp}{\hspace{\fill} $\Box$ \par \bigskip}

\begin{document}
\centerline{\LARGE A Clean Approach to Rational Cubic Residues}

\bigskip\centerline{\Large Sam Vandervelde}

\vspace{5mm} \begin{quote}
\centerline{\textbf{Abstract}}
\noindent In 1958 E.\ Lehmer found an explicit description of those primes $p$ for which a given prime $q$ is a cubic residue.  Her result essentially states that if one writes $4p=L^2+27M^2$, then $q$ is a cubic residue if and only if \[ \frac{M}{L} \equiv \frac{t^2-1}{t^3-9t} \bmod q \] for some integer~$t$.  In this paper we demonstrate that a similar result may be obtained for cubic nonresidues, yielding a cubic character for fixed $p$ that provides an effective means for ascertaining whether or not an arbitrary integer $c$ is a cubic residue modulo~$p$.  As an illustration of this technique, we determine whether 1982 is a cubic residue modulo the 131-digit prime $p=\frac{1}{4}(3^{19}+5^{82})$, a question which, for all practical purposes, is impossible to answer with Lehmer's original criterion.
\end{quote}

\section{Motivation}
Let $m\ge 2$ be a positive integer.  We say that an integer $c$ is a cubic residue modulo $m$ if there is a solution to the congruence $x^3\equiv c \bmod m$.  The study of cubic residues dates back centuries, at least to the mid 1700's when Euler conjectured that 2 is a cubic residue modulo a prime $p\equiv 1\bmod 3$ if and only if $p$ can be written in the form $A^2+27B^2$.  Gauss later proved this statement as a consequence of his theory of cubic reciprocity, which he never published.  It was left to Eisenstein and Jacobi to provide a documented development of this theory, from which results on cubic residues follow.  These historical tidbits, along with a great many more, appear in Cox's book \cite{Cox}.

It is straightforward to show that $x^3\equiv c \bmod m$ has a solution if and only if $x^3\equiv c \bmod p$ has a solution for each prime factor $p$ of~$m$.  (With one exception --- one must consider $x^3\equiv c \bmod 9$ when $m$ is divisible by~$3^k$ with $k>1$, which is solvable exactly when $c\equiv \pm 1\bmod 9$.)  Since every integer $c$ is a cubic residue modulo $p$ for any prime $p\equiv 2\bmod 3$ and for $p=3$, the question is reduced to studying the congruence $x^3\equiv c \bmod p$ for primes $p\equiv 1\bmod 3$.

In 1920 Cunningham and Gosset~\cite{CunGos} compiled tables to settle this question for all integers~$c$ whose prime factors did not exceed 50.  Using a formulation of cubic reciprocity due to P\'epin~\cite{Pep} (which bears a surprisingly close resemblance to more modern notation), they demonstrated that for a given prime~$q$, the value of a certain cubic residue symbol depended only on the ratio $\frac{M}{L}$ modulo~$q$, where $L$~and~$M$ are integers for which $4p=L^2+27M^2$.  The values corresponding to the various possible ratios were obtained by factoring polynomials in $L$~and~$M$ whose degree was proportional to~$q$, making the process increasingly impractical as $q$ grew larger.

This difficulty was circumvented by Lehmer in~\cite{Leh}, where she developed a different sort of criterion that predicted for which primes~$p$ a given prime~$q$ is a cubic residue.  We present a slightly more tidy version of Lehmer's result here; its derivation from the original appears in a later section.

\begin{thm} Let $p\equiv 1\bmod 3$ be a prime and write $4p=L^2+27M^2$.  Then a given prime $q\ne p$ is a cubic residue modulo $p$ if and only if
\begin{equation} (t-1)(t+1)L \equiv t(t-3)(t+3)M \bmod q \label{lehmer} \end{equation}
for some integer~$t$. \label{base}\end{thm}

The result implies that the primes~$p$ for which $q$ is a cubic residue lie on certain lines through the origin modulo $q$ (in the {$(L,M)$-plane) whose slopes are given by
\[ \frac{M}{L} \equiv \frac{t^2-1}{t^3-9t} \bmod q \]
as $t$ runs through all distinct values modulo~$q$.  The values $t=\pm 1$ give horizontal lines, while $t=0$ and $\pm 3$ correspond to vertical lines.  For $q=2$, 3, 5, and 7 these are the only lines that occur; in other words, in these cases $q$~is a cubic residue modulo $p$ precisely when $4p=L^2+27M^2$ has either $L$ or $M$ divisible by~$q$.  (For example, when $q=2$ one can now recover Euler's conjecture mentioned above.)  These lines remain for $q\ge 11$, but oblique lines are also present beyond this point.  It should be noted that for $p\equiv 1\bmod 3$ there are exactly four ways to write $4p=L^2+27M^2$, which differ only in sign.  The theorem applies equally well to any representation, as can be seen by replacing $t$ by $-t$ in~(\ref{lehmer}).

It is interesting to observe that for $q\equiv 2 \bmod 3$, exactly one-third of the $q+1$ possible slopes (including vertical lines) produce primes for which $q$ is a cubic residue, because each such slope is produced by precisely three distinct values of~$t$.  (One must include the value $t=\infty$, yielding a slope of~0, for this to work out properly.)  The same is true, with a few qualifications, when $q\equiv 1\bmod 3$.  We will investigate these considerations more closely in a later section.

It is natural at this point to ask whether it is possible to obtain a complementary statement for cubic non-residues in order to extend the results of Cunningham and Gosset.  Happily, the answer is in the affirmative.  The purpose of the present paper is to provide a clean, self-contained description of how such a result may be formulated, with the goal of ultimately presenting a ``rational'' version of our theory that may be stated without reference to any rings other than the ordinary integers.

Finally, we must acknowledge the progress made by many others on cubic residues and congruences, especially the more recent work of Lehmer, Williams, Spearman, and Sun. (See \cite{Leh2}, \cite{Leh3}, \cite{SpeWil}, \cite{Sun}, and \cite{Wil}.)  The paper of Sun is of particular interest, since it independently develops  some of the same themes, albeit with different emphases and purposes, as this article.

\section{Preliminaries}
We will require the machinery of cubic reciprocity, so we briefly recall some of the relevant definitions and results, as presented in Ireland and Rosen~\cite{IreRos}.  Let $R$ refer to the unique factorization domain $\mathbb{Z}[\omega]=\{a+b\omega\mid a,b \in \mathbb{Z}\}$ with units $\pm1$, $\pm\omega$, and $\pm\omega^2$, where $\omega=e^{2\pi i/3}$.  In what follows, $q$ will always denote a rational prime, i.e. a positive integer that is prime within~$\mathbb{Z}$.  We know that if $q\equiv 2\bmod 3$, it is also prime within~$R$, while if $q\equiv 1\bmod 3$, then $q$ splits into a product $q=\rho\bar{\rho}$ for primes $\rho,\bar{\rho}\in R$.  We also have $3=-\omega^2(1-\omega)^2$ for the prime $1-\omega$.  Furthermore, $p$~will always refer to a rational prime satisfying $p\equiv 1\bmod 3$, so that we can write $p=\pi\bar{\pi}$ for primes $\pi,\bar{\pi}\in R$.  

We say that a prime $\pi=a+b\omega\in R$ is primary if $b$ is divisible by 3 and complex if $b\ne 0$.  (Our terminology differs slightly from~\cite{IreRos} here.)  When $p\equiv 1\bmod 3$, there are exactly four primary primes in $R$ that divide~$p$, all complex.  If $\pi$ is one such prime then the others are $-\pi$ and~$\pm\bar{\pi}$, and $p=\pi\bar{\pi}=(-\pi)(-\bar{\pi})$.  We say that these primes arise from~$p$.  (They are the primary primes with norm~$p$.)  Observe that in this case
\[ 4p = 4\pi\bar{\pi} = 4a^2-4ab+4b^2 = (2a-b)^2+27({\textstyle\frac{b}{3}})^2, \]
so letting $L=2a-b$ and $M=\frac{b}{3}$ yields the four decompositions of $4p$ used in~(\ref{lehmer}), one for each primary divisor $\pi$ of~$p$.  On the other hand, if $q\equiv 2\bmod 3$, then $\pm q$ are the two primary primes in~$R$ dividing (arising from)~$q$.  The only primes in $R$ not yet considered are the divisors of 3, such as $1-\omega$, but none of these are primary.

For each primary prime $\pi\in R$ one can define~$\chi_{\pi}$, a \textit{cubic residue character} modulo $\pi$ on~$R$, as follows.  First, if $\pi$ divides $\alpha\in R$, we set $\chi_{\pi}(\alpha)=0$.  If $\pi$ is complex, then for any $\alpha\in R$ not divisible by~$\pi$ the power $\alpha^{(p-1)/3}$ will be congruent to exactly one of 1, $\omega$, or $\bar{\omega}$ modulo~$\pi$, where $p=\pi\bar{\pi}$ as usual.  One defines $\chi_{\pi}(\alpha)$ to be this cube root of unity.  Otherwise $\pi$ is of the form~$\pm q$, where $q\equiv 2\bmod 3$ is a rational prime.  In this case we use the power $\alpha^{(q^2-1)/3}$ to define $\chi_{\pi}(\alpha)$ in the same manner as before.  Note that in either case $\chi_{\pi}(1)=
\chi_{\pi}(-1)=1$.

It will be helpful to recapitulate several properties of these characters for later use.  All but the last follow more or less directly from the definition.

\medskip
\textit{i.} $\chi_{\pi}(\alpha\alpha') = \chi_{\pi}(\alpha)\chi_{\pi}(\alpha')$.

\textit{ii.} If $\alpha\equiv \alpha'\bmod \pi$, then $\chi_{\pi}(\alpha)=\chi_{\pi}(\alpha')$.

\textit{iii.} $\chi_{\pi}(\alpha)=1$ if and only if $\alpha$ is a cubic residue modulo~$\pi$.

\textit{iv.} The complex conjugate of $\chi_{\pi}(\alpha)$ is $\chi_{\bar{\pi}}(\bar{\alpha})$.

\textit{v.} If $q\equiv 2\bmod 3$ and $k\in\mathbb{Z}$ is not divisible by~$q$, then $\chi_q(k)=1$.

\textit{vi.} If $\pi=a+b\omega$ then $\chi_{\pi}(1-\omega)=\omega^{2n}$, where $a=\pm(3n-1)$.

\medskip
With this background we can now state the law of cubic reciprocity.

\begin{thm} If $\pi$ and $\rho$ are primary, arising from different rational primes, then
\begin{equation} \chi_{\pi}(\rho) = \chi_{\rho}(\pi). \label{cr} \end{equation} \end{thm}
It is arguably the most elegant of all the reciprocity laws.

\section{Constructing the Tables}
Given a rational prime~$q$, our first task is to address the ordered pairs $(L,M)$ corresponding to primes $p$ for which $q$ is not a cubic residue.  Loosely speaking, we wish to assign to each such pair a value, either $\omega$ or~$\bar{\omega}$, that indicates which type of non-residue $q$ should be modulo~$p$.  We have seen that the four pairs of integers $(L,M)$ for which $4p=L^2+27M^2$ correspond to the four pairs $(a,b)$ for which $\pi=a+b\omega$ is a primary divisor of $p$ via $L=2a-b$ and $M=\frac{b}{3}$.  Therefore we define 
\begin{equation} F_q(L,M) = \chi_{\pi}(q). \label{define} \end{equation}
For example, it is a simple exercise to verify that $F_2(1,1)=F_2(-1,-1)=\omega$ while $F_2(1,-1)=
F_2(-1,1)=\bar{\omega}$.  This may seem inconsistent at first, since all four $(L,M)$ pairs correspond to the same prime $p=7$.  However, we are not defining a single character for $p=7$ at this point; rather, we are employing four separate characters, one for each primary divisor of 7.

We know that $\chi_{\pi}(q)=1$ means that $q$ is a cubic residue modulo~$\pi$.
Thus $F_q(L,M)=1$ implies that $x^3\equiv q \bmod \pi$ for some $x\in R$.  Since the set $\{0,1,\ldots, 
p-1\}$ represents a complete set of residues modulo~$\pi$, we may assume that $x$ is an integer.  Conjugating gives $x^3\equiv q \bmod \bar{\pi}$, hence $x^3\equiv q \bmod p$.  The converse is clear, so we have shown that $F_q(L,M)=1$ if and only if $q$ is a cubic residue modulo~$p$.  Lehmer's theorem provides an alternate way of determining those pairs $(L,M)$ for which $F_q(L,M)=1$ without computing the cubic characters $\chi_{\pi}(q)$.  Our goal is to provide a similar description of those pairs for which $F_q(L,M)=\omega$ or~$\bar{\omega}$.

\begin{prop} Let $q\ge 5$ be prime.  Then the value of $F_q(L,M)$ is constant along lines through the origin, modulo~$q$.  Furthermore, pairs of lines whose slopes are negatives of one another have conjugate values; that is, $F_q(L,-M)=\overline{F_q(L,M)}$. \label{lines}\end{prop}

Before beginning the proof, we remark that there is a potential difficulty involving points $(L,M)$ congruent to the origin modulo~$q$, since all the lines described above pass through such a point.  However, this would mean that $q|L$ and $q|M$, in which case $L^2+27M^2$ could not equal four times a prime.  In other words, $F_q(L,M)$ is never defined at such points.  The one exception occurs
for $q=2$, since $4\cdot31=4^2+27\cdot2^2$, for instance.  It turns out that one must examine the values of $F_2(L,M)\bmod 4$ rather than mod 2.

\medskip
\pf The points $(L,M)$ on a particular line through the origin, modulo~$q$, are precisely those pairs satisfying $L\equiv kl \bmod q$ and $M\equiv km \bmod q$, where $l$ and $m$ are fixed integers (not both divisible by $q$) and $k$ takes on the values $k=0$, 1, \ldots,~$q-1$.  Taking $l=1$, for instance, yields a line of slope~$m$, while choosing $l=0$ and $m=1$ gives a vertical line.  Since 
\[ a={\textstyle\frac{1}{2}}(L+3M)=k\cdot{\textstyle\frac{1}{2}}(l+3m), \qquad b=3M=k\cdot 3m, \]
the corresponding pairs $(a,b)$ also lie on a certain line, the one for which $a\equiv ka'\bmod q$ and $b\equiv kb'\bmod q$, where $a'=\frac{1}{2}(l+3m)$ and $b'=3m$.  Recall that $\pi=a+b\omega$.  We set $\pi'=a'+b'\omega$ so that $\pi\equiv k\pi'\bmod q$.

To apply cubic reciprocity, we must consider $q\equiv 1, 2\bmod 3$ separately.  In the latter case $q$ is primary, so $F_q(L,M) = \chi_{\pi}(q) = \chi_q(\pi)$.  Using the fact that $\chi_q$ is multiplicative, we deduce that $F_q(L,M) = \chi_{q}(k) \chi_q(\pi')$.  We can disregard $k\equiv 0\bmod q$, since $q$ would divide both $L$~and~$M$.  Hence $\chi_q(k)=1$, which means that $F_q(L,M)=\chi_q(\pi')$ is constant.  The argument is nearly identical when $q\equiv 1\bmod 3$.  In this case we first write $q=\rho\bar{\rho}$ where $\rho$ and $\bar{\rho}$ are primary.  Therefore 
\[ F_q(L,M) = \chi_{\pi}(q) = \chi_{\pi}(\rho)\chi_{\pi}(\bar{\rho}) = \chi_{\rho}(\pi) \chi_{\bar{\rho}}(\pi). \]
Since $\pi\equiv k\pi'\bmod q$ we know the same is true modulo $\rho$ or $\bar{\rho}$.  Hence
\[ F_q(L,M) = \chi_{\rho}(k) \chi_{\bar{\rho}}(k) \chi_{\rho}(\pi') \chi_{\bar{\rho}}(\pi'). \]
But the first two factors are conjugates, so their product is 1, while the second two depend only on~$\pi'$.  We again conclude that $F_q(L,M)$ is constant along the particular line under consideration.

Finally, note that if the point $(L,M)$ corresponds to a prime~$\pi$, then $(L,-M)$ corresponds to~$\bar{\pi}$, since $\bar{\pi}=a+b\bar{\omega}=(a-b)-b\omega$ is associated with the point 
$(2(a-b)-(-b),\frac{-b}{3}) = (2a-b,-\frac{b}{3})$, which is $(L,-M)$.  Therefore 
\[ F_q(L,-M)=\chi_{\bar{\pi}}(q)= \overline{\chi_{\pi}(q)}=\overline{F_q(L,M)},\] as desired. \fp

Similar results hold for $q=2$ and $q=3$, but the previous approach cannot be employed.  (The expression $\frac{1}{2}(L+3M)$ makes no sense mod 2, and when $q=3$ cubic reciprocity does not apply because 3 has no primary divisors.)  Therefore we modify the argument for $q=2$, and simply
compute $F_3(L,M)$ directly.

\begin{prop} The values of $F_2(L,M)$ and $F_3(L,M)$ are given by the following criteria:
\[ F_2(L,M) = \left\{ \begin{array}{cl} 1 & L, M\equiv 0 \bmod 2 \\
\omega & L\equiv M \bmod 4 \\ 
\bar{\omega} & L\equiv -M \bmod 4 \end{array} \right. \quad
F_3(L,M) = \left\{ \begin{array}{cl} 1 & M\equiv 0 \bmod 3 \\
\omega & L\equiv -M \bmod 3 \\ 
\bar{\omega} & L\equiv M \bmod 3 \end{array} \right. .\] \label{twothree}\end{prop}

\pf Since 2 is a primary prime in~$R$ we may employ cubic reciprocity to write $F_2(L,M)=\chi_{\pi}(2)=\chi_2(\pi)$.  By definition, $\chi_2(\pi)\equiv \pi\bmod 2$.  Writing $\pi=a+b\omega$, we find that $\chi_2(\pi)=\bar{\omega}=-1-\omega$ if and only if $a\equiv b\equiv -1\bmod 2$.  Since $L=2a-b$, this translates to $L\equiv 1, 3\bmod 4$, depending upon whether $b\equiv 1, 3\bmod 4$, respectively.  But $M=\frac{b}{3}\equiv -b\bmod 4$, so $M\equiv 3, 1\bmod 4$ according as to $b\equiv 1, 3\bmod 4$.  Regardless, $L\equiv -M\bmod 4$, as stated.  The analysis of $\chi_2(\pi)=1$ or $\omega$ is also straightforward, as the reader may verify.

To compute $F_3(L,M)$ we use the fact that $3=-\omega^2(1-\omega)^2$ to write
\[ F_3(L,M) = \chi_{\pi}(3) = \chi_{\pi}(-1)\cdot [\chi_{\pi}(\omega)]^2\cdot [\chi_{\pi}(1-\omega)]^2. \]
By definition, $\chi_{\pi}(\omega)=\omega^{(p-1)/3}$.  Using the properties of $\chi_{\pi}$ presented earlier, we find that
\[ F_3(L,M) = \omega^{(2p-2+3n)/3}, \]
where $p=\pi\bar{\pi}=a^2-ab+b^2$ and $a=\pm(3n-1)$.  Hence we must determine the exponent mod~3, which means we must evaluate $2p-2+3n\bmod 9$.  For example, suppose that $L\equiv M\equiv 1\bmod 3$.  Then $b=3M\equiv 3\bmod 9$, while $a=\frac{1}{2}(L+3M)\equiv 2\bmod 3$.  First note that $p\equiv a^2+3\bmod 9$.  In addition, $a$ is of the form $3n-1$ (as opposed to $-(3n-1)$), so $3n=a+1$.  Therefore we may write
\[ 2p-2+3n \equiv (2a^2+6)-2+a+1 \equiv (2a-1)(a+1)+6 \equiv 6 \bmod 9, \]
since both $(a+1)$ and $(2a-1)$ are divisible by 3.  We conclude that $F_3(L,M)=\omega^2=\bar{\omega}$, as claimed.  The remaining possibilities may be ascertained in an analogous manner, leading in each case to the result stated above.  \fp

In light of the periodicity of the values of $F_q(L,M)$ for a given prime~$q$, it makes sense to define a function $f_q:\mathbb{Z}/q\mathbb{Z}\times\mathbb{Z}/q\mathbb{Z}\rightarrow \{1,\omega,\bar{\omega}\}$ by setting $f_q(l,m)=F_q(L,M)$, for $l\equiv L\bmod q$ and $m\equiv M\bmod q$, whenever $F_q(L,M)$ is defined and non-zero.  (Remembering, of course, to operate mod~4 when $q=2$.)  This function is well-defined by the preceding propositions and encapsulates all the information about~$F_q$.  Although we will not need it, one can show that if $f_q(l,m)$ is defined for one point on a line through the origin, then it will be defined for every point on that line aside from the origin.

It is informative to present $f_q$ as a $q\times q$ (or $4\times 4$) table of values, which we have done for $q=2$, 3, 5, and 7 in Figure~\ref{tables}.  A~`$*$' appears in the table to indicate that $F_q(L,M)$ is never defined (or possibly equals 0) for that particular point, because the expression $L^2+27M^2$ is never equal to four times a prime (other than possibly $q$) for $l\equiv L\bmod q$ and $m\equiv M\bmod q$.  As expected, $f_q(0,0)=*$, but there are other occurrences as well.  Thus $L^2+27M^2$ is odd when $L$ and $M$ have opposite parity, is divisible by 16 when $L\equiv M\equiv 0, 2\bmod 4$, is divisible by 9 when $L\equiv 0\bmod 3$, and is divisible by 7 when $L\equiv \pm M\bmod 7$.  The latter instance implies that $\frac{1}{4}(L^2+27M^2)$ is composite unless it equals 7, so that $F_7(L,M)$ is not defined in this case except for $F_7(\pm 1,\pm 1)=0$.

For small odd values of $q$ the tables are relatively easy to construct.  There are $q+1$ lines through the origin, on which $f_q$ is constant.  The horizontal and vertical line contain 1's by Lehmer's theorem and pairs of lines with negative slopes contain conjugate values by Proposition~\ref{lines}.  Hence to complete the table for $f_5$ we need only determine its value at two strategic points.  For instance, computing $F_5(1,1)=\chi_{2+3\omega} (5)=\omega$ and $F_5(4,2)=\chi_{5+6\omega}(5)=\bar{\omega}$ will suffice.  For the same reasons we need only compute two values of $F_7$ to deduce the entire table.  Notice that in each of the four tables exactly one-third of the available lines, and hence one-third of the non-starred entries, are filled with each possible value, either 1, $\omega$, or~$\bar{\omega}$.

\begin{figure}\renewcommand{\arraystretch}{1.2}
\centerline{\begin{tabular}{c}
\begin{tabular}{r|c|c|c|} \multicolumn{1}{r}{} & 
\multicolumn{3}{c}{$f_3(l,m)$} \\ \cline{2-4}
2 & $*$ & $\omega$ & $\bar{\omega}$ \\ \cline{2-4} 
1 & $*$ & $\bar{\omega}$ & $\omega$ \\ \cline{2-4}
0 & $*$ & 1 & 1 \\ \cline{2-4} 
\multicolumn{1}{r}{} & \multicolumn{1}{c}{0} & \multicolumn{1}{c}{1} & 
\multicolumn{1}{c}{2} \end{tabular} \\
\begin{tabular}{r|c|c|c|c|} \multicolumn{1}{r}{} & \multicolumn{4}{c}{$f_2(l,m)$} \\ \cline{2-5} 
3 & $*$ & $\bar{\omega}$ & $*$ & $\omega$ \\ \cline{2-5}
2 & 1 & $*$ & $*$ & $*$ \\ \cline{2-5} 
1 & $*$ & $\omega$ & $*$ & $\bar{\omega}$ \\ \cline{2-5}
0 & $*$ & $*$ & 1 & $*$ \\ \cline{2-5} 
\multicolumn{1}{r}{} & \multicolumn{1}{c}{0} & \multicolumn{1}{c}{1} & 
\multicolumn{1}{c}{2} & \multicolumn{1}{c}{3} \end{tabular}
\end{tabular}\hspace{1mm}%
\begin{tabular}{r|c|c|c|c|c|} \multicolumn{1}{r}{} & \multicolumn{5}{c}{$f_5(l,m)$} \\ \cline{2-6} 
4 & 1 & $\bar{\omega}$ & $\omega$ & $\bar{\omega}$ & $\omega$ \\ \cline{2-6}
3 & 1 & $\bar{\omega}$ & $\bar{\omega}$ & $\omega$ & $\omega$ \\ \cline{2-6}
2 & 1 & $\omega$ & $\omega$ & $\bar{\omega}$ & $\bar{\omega}$ \\ \cline{2-6} 
1 & 1 & $\omega$ & $\bar{\omega}$ & $\omega$ & $\bar{\omega}$ \\ \cline{2-6}
0 & $*$ & 1 & 1 & 1 & 1 \\ \cline{2-6} 
\multicolumn{1}{r}{} & \multicolumn{1}{c}{0} & \multicolumn{1}{c}{1} & 
\multicolumn{1}{c}{2} & \multicolumn{1}{c}{3} & \multicolumn{1}{c}{4} \end{tabular}
\hspace{4mm}%
\begin{tabular}{r|c|c|c|c|c|c|c|} \multicolumn{1}{r}{} & \multicolumn{7}{c}{$f_7(l,m)$} \\ \cline{2-8} 
6 & 1 & $*$ & $\bar{\omega}$ & $\omega$ & $\bar{\omega}$ & $\omega$ & $*$ \\ \cline{2-8}
5 & 1 & $\bar{\omega}$ & $*$ & $\omega$ & $\bar{\omega}$ & $*$ & $\omega$ \\ \cline{2-8}
4 & 1 & $\omega$ & $\omega$ & $*$ & $*$ & $\bar{\omega}$ & $\bar{\omega}$ \\ \cline{2-8}
3 & 1 & $\bar{\omega}$ & $\bar{\omega}$ & $*$ & $*$ & $\omega$ & $\omega$ \\ \cline{2-8}
2 & 1 & $\omega$ & $*$ & $\bar{\omega}$ & $\omega$ & $*$ & $\bar{\omega}$ \\ \cline{2-8} 
1 & 1 & $*$ & $\omega$ & $\bar{\omega}$ & $\omega$ & $\bar{\omega}$ & $*$ \\ \cline{2-8}
0 & $*$ & 1 & 1 & 1 & 1 & 1 & 1 \\ \cline{2-8} 
\multicolumn{1}{r}{} & \multicolumn{1}{c}{0} & \multicolumn{1}{c}{1} & 
\multicolumn{1}{c}{2} & \multicolumn{1}{c}{3} & \multicolumn{1}{c}{4} 
& \multicolumn{1}{c}{5} & \multicolumn{1}{c}{6} \end{tabular}\hspace*{6mm}}
\renewcommand{\arraystretch}{1}
\caption{The table of values for $f_q(l,m)$ when $q=2$, 3, 5, and 7.}
\label{tables}\end{figure}

As a simple application, let us determine whether or not $490$ is a cubic residue mod~63601.  First, one finds that $4\cdot 63601=19^2+27\cdot 97^2$.  Hence we may use any of the four pairs $(\pm 19,\pm 97)$ in our computations, as long as we use the same pair consistently.  Taking $(L,M)=(19,97)$ we have $F_2(19,97)=f_2(3,1)=\bar{\omega}$, $F_5(19,97)=f_5(4,2)=\bar{\omega}$, and $F_7(19,97)=f_7(5,6)=\omega$.  The product $\bar{\omega}\cdot\bar{\omega}\cdot\omega^2$ corresponding to $490=2\cdot 5\cdot 7^2$ is equal to 1, so we conclude that 490 is a cubic residue.  This process works because we are actually computing $\chi_{155+291\omega}(490)$ and using the fact that 490 is a cubic residue mod~63601 if and only if the value of this character equals 1.

\section{Determining the $\omega$-slopes}
Given a prime~$q$, we now ascertain the pairs $(L,M)$ for which $F_q(L,M)=\omega$ by describing the slopes of the lines containing these points.  As in Theorem~\ref{base}, this is accomplished with a rational function of degree three, but unlike that result, no single rational function will suffice for all primes.  We begin by describing the $\omega$-slopes when 
$q\not\equiv\pm1\bmod 9$.  The following lemma will be useful.

\begin{lem} Let $q\ne 3$ be a rational prime relatively prime to a given number $a+b\omega\in R$.  Then there exists $\lambda\in R$ primary with $\lambda\equiv a+b\omega \bmod q$. \end{lem}

\pf First choose $a'\equiv a\bmod q$ and $b'\equiv b\bmod q$ so that $b'$ is divisible by 3 but $a'$ is not.  Then consider the sequence $\{a'+b'\omega + k(3q)\mid k\in\mathbb{Z}\}$.  Clearly $3q$ is relatively prime to $a'+b'\omega$, so by Dirichlet's theorem on primes in arithmetic progressions (for number fields) some element $\lambda$ of the sequence is prime.  By construction $\lambda$ is primary and in addition $\lambda\equiv a+b\omega\bmod q$, as desired. \fp 

\begin{prop} Let $q>3$ be a rational prime satisfying $q\equiv\pm 2\bmod 9$, and write $4p=L^2+27M^2$ for a prime $p\equiv 1\bmod 3$, $p\ne q$.  Then $F_q(L,M)=\omega$ if and only if
\begin{equation} (t^3-3t^2-9t+3) L \equiv -3(t^3+9t^2-9t-9) M \bmod q \label{omeqn} \end{equation}
for some integer~$t$.  When $q\equiv\pm 4\bmod 9$ the same result holds without the negative sign. \label{omega}\end{prop}

\pf Let $\pi=a+b\omega$ where $L=2a-b$, $M=\frac{b}{3}$, so that $p=\pi\bar{\pi}$ and $\pi$~is primary, as usual.  We wish to find conditions on $L$~and~$M$ equivalent to $F_q(L,M)=\omega$.  This means that $\chi_{\pi}(q)=\omega$, by definition.  If $q\equiv 2\bmod 9$, so that $q$ is also primary, then $\chi_q(\pi)=\omega$ by cubic reciprocity.  However, we know that $\chi_q(\bar{\omega})= \bar{\omega}^{(q^2-1)/3} =\bar{\omega}$, so that
\[ 1=\omega\bar{\omega}=\chi_q(\pi)\chi_q(\bar{\omega})=\chi_q(\bar{\omega}\pi). \]
At this point we would like to invoke cubic reciprocity again in order to apply Lehmer's criterion, but 
$\bar{\omega}\pi=(b-a)-a\omega$ is no longer primary.  However, there is a primary prime $\lambda\equiv\bar{\omega}\pi\bmod q$ by the lemma; hence $\chi_q(\lambda)=\chi_{\lambda}(q)=1$.  If $(L',M')$ is the pair corresponding to~$\lambda$, then 
\[ L'\equiv 2(b-a)-(-a)\equiv 2b-a\bmod q,\quad M'\equiv {\textstyle\frac{1}{3}}(-a)\bmod q. \]
One confirms that $L'\equiv -\frac{1}{2}L+\frac{9}{2}M\bmod q$ while $M'\equiv -\frac{1}{6}L-\frac{1}{2}M\bmod q$.  According to Theorem~\ref{base}, $\chi_{\lambda}(q)=1$ occurs if and only if
\[ (t-1)(t+1)L' \equiv t(t-3)(t+3)M' \bmod q \] 
for some integer~$t$, which reduces to equation~(\ref{omeqn}) upon writing $L'$~and~$M'$ in terms of $L$~and~$M$. 

If $q\equiv -2\bmod 9$, then $q$ is not prime in $R$ and we must write $q=\rho\bar{\rho}$ for a primary prime $\rho\in R$.  In this case $F_q(L,M)=\chi_{\pi}(q)=\omega$ becomes 
\[ \omega=\chi_{\pi}(\rho)\chi_{\pi}(\bar{\rho})=\chi_{\rho}(\pi)\chi_{\bar{\rho}}(\pi). \]
We find that $\chi_{\rho}(\bar{\omega})=\chi_{\bar{\rho}}(\bar{\omega})=\bar{\omega}^{(q-1)/3}=\omega$ this time; multiplying these equalities into the above equation yields $\chi_{\rho}(\bar{\omega}\pi) \chi_{\bar{\rho}}(\bar{\omega}\pi) = 1$.  Choosing $\lambda$ primary with $\lambda\equiv \bar{\omega}\pi\bmod q$ as before then leads to $\chi_{\rho}(\lambda) \chi_{\bar{\rho}}(\lambda) = 1$, or $\chi_\lambda(q)=1$ upon using cubic reciprocity.  One then proceeds as before to obtain equation~(\ref{omeqn}).

The analysis of the case $q\equiv -4\bmod 9$ is nearly identical to $q\equiv 2\bmod 9$, except that $\chi_q(\omega)=\omega$ this time, so that one considers $\omega\pi=-b+(a-b)\omega$ rather than $\bar{\omega}\pi$.  Now $L'=-a-b=-\frac{1}{2}L-\frac{9}{2}M$ and $M'=\frac{1}{3}(a-b)=\frac{1}{6}L -\frac{1}{2}M$, which leads to the congruence
\[ (t^3+3t^2-9t-3) L \equiv 3(t^3-9t^2-9t+9) M \bmod q. \]
Replacing $t$ by $-t$ and negating both sides then produces~(\ref{omeqn}) sans negative sign, as claimed.  Finally, combining these steps with previous ideas gives the same result when $q\equiv 4\bmod 9$.  It is relatively straightforward to reverse this argument to establish the converse.  However, we omit the details in favor of a neater approach, outlined below.  This completes the proof.  \fp

The foregoing proposition identifies the slopes of the lines through the origin on which $F_q$ (or $f_q$) equals either $\omega$ or~$\bar{\omega}$, when $q\not\equiv\pm 1\bmod 9$.  Thus when $q\equiv\pm 2\bmod 9$ the lines with slope
\begin{equation} \frac{M}{L} \equiv -\frac{t^3-3t^2-9t+3}{3(t^3+9t^2-9t-9)} \bmod q \label{slope}\end{equation}
pass through points for which $F_q(L,M)=\omega$; while the lines whose slopes are the negatives of these satisfy $F_q(L,M)=\bar{\omega}$.  The situation is reversed for $q\equiv\pm 4\bmod 9$.  It is interesting to note that when $q\equiv\pm 1\bmod 9$ these lines predict the primes for which $q$ is a cubic residue; in other words, the rational function above duplicates the slopes given by Lehmer's theorem.  This is to be expected in light of the above proof and the fact that $\chi_q(\omega)=1$ (or $\chi_{\rho}(\omega)=1$) in this case.

It is instructive to confirm Proposition~\ref{omega} for $q=5$ and $q=7$, since we already have tables of values for $f_q(l,m)$.  Letting $t=0$, 1, 2, 3, 4, and $\infty$ in~(\ref{slope}) mod~5 yields slopes of 4, 3, 4, 4, 3, and~3, respectively.  These are the $\bar{\omega}$-slopes for~$f_5$, as predicted by the proposition.  On the other hand, taking $t=0$, 1, \ldots, 6, and~$\infty$ mod~7 gives slopes of 4, 2, 6, 4, 4, 1, 2, and~2, respectively.  These are the $\omega$-slopes for~$f_7$, as claimed, along with the lines on which $f_7$ is not defined, i.e.\ where the proposition does not apply.  Observe that each $\omega$ or $\bar{\omega}$-slope occurs exactly three times.  This behavior is typical of the rational functions that compute such slopes, which we now introduce.

Let $\gamma\in R$ be written as $\gamma=c+d\omega$.  To ease notation, we also set $C=6c-3d$ and $D=d$.  We then define
\begin{equation} g_{\gamma}(t) = -\frac{Dt^3-Ct^2-9Dt+C}{Ct^3+27Dt^2-9Ct-27D}. \label{defn}
\end{equation}
This family of functions possesses several interesting properties.

\begin{prop} The functions $g_{\gamma}$ satisfy $g_{-\gamma}(t)= g_{\gamma}(t)$ and $g_{\bar{\gamma}}(t)=-g_{\gamma}(-t)$.  If $q\equiv 2\bmod 3$ is a rational prime, then $g_{\gamma}(t)$ assumes exactly $\frac{1}{3}(q+1)$ distinct values when $t=0$, $1$, \ldots, $q-1$, and $\infty$; each is attained by three different values of~$t$.  Similarly, if $q\equiv 1\bmod 3$ then $g_{\gamma}(t)$ assumes exactly $\frac{1}{3}(q-1)+2$ distinct values; $\frac{1}{3}(q-1)$ of them each coming from three different values of~$t$, along with the two values $\pm\frac{1}{3\sqrt{-3}}$, each attained once. \label{range} \end{prop}

\pf It is clear that $g_{-\gamma}(t)=g_{\gamma}(t)$.  Next observe that $C=6\,\mathrm{Re}(\gamma)$ while $D=\frac{2}{\sqrt{3}}\,\mathrm{Im} (\gamma)$, so that conjugating $\gamma$ fixes $C$ and negates~$D$.  This explains the fact that $g_{\bar{\gamma}}=-g_{\gamma}(-t)$.  To analyze the range of $g_{\gamma}(t)$ modulo a prime~$q$, we write
\begin{equation} g_{\gamma}(t) = \frac{Cg(t)-D}{27Dg(t)+C}, \quad g(t)=g_1(t)=\frac{t^2-1}{t^3-9t},
\label{altform} \end{equation}
where $g(t)$ is the function appearing in Lehmer's theorem.  Now the transformation $h:t\mapsto\frac{t-3}{t+1}$ has order three, since applying it repeatedly yields
\[ t\mapsto\frac{t-3}{t+1}\mapsto\frac{t+3}{-t+1}\mapsto t. \]
But we have
\[ g(t) = \left[ -(t)\left(\frac{t-3}{t+1}\right)\left(\frac{t+3}{-t+1}\right) \right]^{-1}, \]
so $g(t)$ is invariant under~$h$.  Hence $g_{\gamma}(t)$ is also, because of~(\ref{altform}).  Thus $g_{\gamma}(t)$ assumes the same value modulo~$q$ at~$t$, $\frac{t-3}{t+1}$, and $\frac{t+3}{-t+1}$.  When $q\equiv 2\bmod 3$ these three numbers are distinct modulo~$q$, since the congruence of any two is equivalent to $t^2\equiv -3\bmod q$, which does not occur.  Since $g_{\gamma}(t)$ may assume the same value for at most three distinct values of~$t$, we deduce that $g_{\gamma}(t)$ takes on exactly $\frac{1}{3}(q+1)$ different values in the manner claimed.

On the other hand, when $q\equiv 1\bmod 3$ then $t=\frac{t-3}{t+1}=\frac{t+3}{-t+1}$ twice, for the two values of~$t$ satisfying $t^2\equiv -3\bmod q$.  (Otherwise $g_{\gamma}(t)$ is three-to-one as before.)  Let $\sqrt{-3}$ refer to one of these values modulo~$q$.  We find that
\[ g_{\gamma}(\sqrt{-3}) \equiv \frac{4C-12D\sqrt{-3}}{12C\sqrt{-3}+108D} \equiv \frac{1}{3\sqrt{-3}} \bmod q, \]
upon cancelling a common factor of $C\sqrt{-3}+9D$.  In the same manner, $g_{\gamma}(-\sqrt{-3}) \equiv -\frac{1}{3\sqrt{-3}} \bmod q$, as asserted. \fp

By the proposition, the values of $g_{\bar{\gamma}}(t)$ are the negatives of the values taken by~$g_{\gamma}(t)$.  In addition, the exceptional values of $g_{\gamma}(t)$ are exactly the slopes for which $F_q(L,M)$ is not defined, since
\[ L^2+27M^2\equiv 0\bmod q \quad\iff\quad \frac{M}{L} \equiv \pm\frac{1}{3\sqrt{-3}} \bmod q. \]
We are now in a position to provide a neat conclusion to Proposition~\ref{omega}.  Observe that
$g_{\omega}(t)$ and $g_{\bar{\omega}}(t)$ are the rational functions giving the $\omega$ and $\bar{\omega}$-slopes that were encountered there.  When $q\equiv 2\bmod 3$, we have seen that each of~$g(t)$, $g_{\omega}(t)$, and $g_{\bar{\omega}}(t)$ take on exactly $\frac{1}{3}(q+1)$ distinct values modulo~$q$.  But these account for all the possible slopes, which proves the only if portion of the statement immediately.  The same logic applies when $q\equiv 1\bmod 3$, once we count the lines on which $F_q(L,M)$ is not defined separately.

The proof of Proposition~\ref{omega} hinged upon the fact that the value of $\chi_q(\omega)$ (or $\chi_{\rho}(\omega)$) could be computed and depended upon $q\bmod 9$.  As long as that value was not equal to 1 we obtained a formula for $\omega$ or $\bar{\omega}$-slopes.  Our next theorem provides an entire family of such formulas, one for each prime $l\equiv1\bmod 3$.  (We have omitted an analysis of $\chi_q(1-\omega)$; it turns out that one obtains the same result as in Proposition~\ref{omega}.)

\begin{thm} Let $l\equiv 1\bmod 3$ be a rational prime, and write $l=\gamma\bar{\gamma}$ where $\gamma, \bar{\gamma}\in R$ are primary.  Then $g_{\gamma}(t)$ gives the slopes, modulo~$q$, on which $F_q(L,M)=1$, $\omega$, or $\bar{\omega}$, and the particular value of $F_q$ occurring along these lines depends upon the value of $q$ modulo~$l$. \label{main} \end{thm}

\pf This argument closely parallels the proof of Proposition~\ref{omega}, so we will be brief.  First write $\gamma=c+d\omega$.  We know that $\chi_{\gamma}(q)\equiv q^{(l-1)/3}\bmod \gamma$, so $\chi_{\gamma}(q)$ depends on the value of $q$ modulo~$\gamma$, and hence modulo~$l$.  Suppose that $\chi_{\gamma}(q)=\bar{\omega}$.  We wish to identify those pairs $(L,M)$ for which $F_q(L,M)=\omega$, where $4p=L^2+27M^2$, $p=\pi\bar{\pi}$, $\pi=a+b\omega$, $L=2a-b$, and $M=\frac{b}{3}$ as usual.  In the case that $q$ is prime in $R$ cubic reciprocity implies that $\chi_q(\pi)=\omega$ and $\chi_q(\gamma)=\bar{\omega}$, hence $\chi_q(\gamma\pi)=1$, or $\chi_{\lambda}(q)=1$ for some primary prime $\lambda\equiv\gamma\pi\bmod q$.  If $(L',M')$ is the pair corresponding to~$\lambda$, then
\[ L'\equiv 2ac-bd-ad-bc \bmod q, \quad M' \equiv {\textstyle\frac{1}{3}}(ad+bc-bd) \bmod q. \]
One then writes $L'$ and $M'$ in terms of $L$ and $M$ to obtain
\[ L' = \left(\frac{2c-d}{2}\right) L - \left(\frac{9d}{2}\right) M, \quad
M' = \left(\frac{d}{6}\right) L + \left(\frac{2c-d}{2}\right) M. \]
Applying Theorem~\ref{base}, we conclude that $(t^2-1)L' \equiv (t^3-9t)M' \bmod q$ for some integer~$t$.  In terms of $L$~and~$M$ this becomes
\[ (Dt^3-Ct^2-9Dt+C)L \equiv -(Ct^3+27Dt^2-9Ct-27D)M \bmod q, \]
where $C=6c-3d$ and $D=d$.  In other words, $F_q(L,M)=\omega$ implies that the slope $\frac{M}{L}$ is of the form~$g_{\gamma}(t)$, modulo~$q$.

Similarly, one finds that $F_q(L,M)=\bar{\omega}$ means that $\frac{M}{L}\equiv g_{\bar{\gamma}}(t) \bmod q$, since $\chi_{\bar{\gamma}}(q)=\omega$.  Of course, $F_q(L,M)=1$ leads to $\frac{M}{L}\equiv g(t) \bmod q$ by Lehmer's result.  Since these three cases account for all of the $q+1$ possible slopes when $q\equiv 2\bmod 3$, by Proposition~\ref{range}, we can make the stronger assertion that $F_q(L,M)=\omega$ if and only if $\frac{M}{L}\equiv g_{\gamma}(t) \bmod q$ for some~$t$.

If $q$ is not prime in $R$ then we write $q=\rho\bar{\rho}$, where $\rho,\bar{\rho}\in R$ are primary.  A few extra steps are required, but just as in Proposition~\ref{omega} one still reaches $\chi_{\lambda}(q)=1$ with $\lambda\equiv\gamma\pi \bmod q$, so $g_{\gamma}(t)$ prescribes the $\omega$-slopes regardless of whether or not $q$ is prime in~$R$.  This concludes the discussion of the case $\chi_{\gamma}(q)=\bar{\omega}$.

On the other hand, if $q$ is a rational prime such that $\chi_{\gamma}(q)=\omega$, then the same reasoning shows that $F_q(L,M)=\omega$ if and only if $\frac{M}{L}\equiv g_{\bar{\gamma}}(t) \bmod q$, while $F_q(L,M)=\bar{\omega}$ is equivalent to $\frac{M}{L}\equiv g_{\gamma}(t) \bmod q$; in other words, the situation is reversed in this case.  Finally, if $q$ satisfies $\chi_{\gamma}(q)= \chi_{\bar{\gamma}}(q)=1$ then Lehmer's Theorem implies that both $g_{\gamma}(t)$ and $g_{\bar{\gamma}}(t)$ give slopes on which $F_q(L,M)=1$, so no formula for $\omega$-slopes is obtained.  This completes the proof. \fp

\section{Rational Cubic Residues}
We now give a pared-down, ``rational'' version of this theory, presented without any reference to the ring $R$ which was instrumental to its development.  As explained at the outset, the question of whether or not the congruence $x^3\equiv c\bmod m$ has a solution may be reduced to determining whether or not $x^3\equiv c\bmod p$ is solvable for a given prime $p\equiv 1\bmod 3$.  This can be done by utilizing a cubic character $\chi$ modulo~$p$, which we now construct.  To begin, set $\chi(0)=\chi(p)=0$ and $\chi(1)=\chi(-1)=1$.  Also let $L$~and~$M$ be the unique \textit{positive\/} integers for which $4p=L^2+27M^2$.  Then given a prime number $q\ne p$, we define $\chi(q)$ according to the following rules.  The first two handle the special cases of $q=2$ and $q=3$.  In the remaining rules we assume that $q\ge 5$.

\bigskip {\setlength{\parindent}{0cm}
(a) If $4|LM$ then $\chi(2)=1$, if $4|(L-M)$ then $\chi(2)=\omega$, and if $4|(L+M)$ then $\chi(2)=\bar{\omega}$.

(b) If $3|M$ then $\chi(3)=1$, if $3|(L+M)$ then $\chi(3)=\omega$, and if $3|(L-M)$ then $\chi(3)=\bar{\omega}$.

(c) If $(t^2-1)L\equiv(t^3-9t)M\bmod q$ for some integer~$t$, then $\chi(q)=1$.

(d) Suppose $q\not\equiv\pm 1\bmod 9$.  Then either $\chi(q)=1$ or
\[ \frac{M}{L} \equiv \pm\frac{t^3-3t^2-9t+3}{3(t^3+9t^2-9t-9)} \bmod q \]
for some integer~$t$.  When $q\equiv\pm 2\bmod 9$ and the plus sign occurs, set $\chi(q)=\bar{\omega}$; with the minus sign define $\chi(q)=\omega$.  If instead $q\equiv\pm 4\bmod 9$ and the plus sign holds, then let $\chi(q)=\omega$, otherwise define $\chi(q)=\bar{\omega}$.

(e) Suppose $q\not\equiv\pm 1\bmod 7$.  Then either $\chi(q)=1$ or
\[ \frac{M}{L} \equiv \pm\frac{t^3-t^2-9t+1}{t^3+27t^2-9t-27} \bmod q \]
for some integer~$t$.  When $q\equiv\pm 2\bmod 7$ and the plus sign occurs, set $\chi(q)=\omega$; with the minus sign define $\chi(q)=\bar{\omega}$.  If instead $q\equiv\pm 3\bmod 7$ and the plus sign holds, then let $\chi(q)=\bar{\omega}$, otherwise define $\chi(q)=\omega$.

(f) Suppose $q\not\equiv\pm 1,\pm 5\bmod 13$.  Then either $\chi(q)=1$ or
\[ \frac{M}{L} \equiv \pm\frac{t^3-5t^2-9t+5}{5t^3+27t^2-45t-27} \bmod q \]
for some integer~$t$.  When $q\equiv\pm 2,\pm 3\bmod 13$ and the plus sign occurs, set $\chi(q)=\omega$; with the minus sign define $\chi(q)=\bar{\omega}$.  If $q\equiv\pm 4, \pm 6\bmod 13$ and the plus sign holds, then let $\chi(q)=\bar{\omega}$, otherwise define $\chi(q)=\omega$.}

\bigskip\noindent Finally, one extends $\chi$ multiplicatively to define $\chi(c)$ when $c$ is not prime.  Thus if $c=\pm q_1^{e_1}\cdots q_n^{e_n}$, then we set $\chi(c)=\chi(q_1)^{e_1}\cdots\chi(q_n)^{e_n}$.

\begin{thm} The construction of $\chi$ just given is self-consistent and defines a cubic character modulo~$p$.  Furthermore, $\chi(c)=1$ if and only if $c$ is a cubic residue modulo~$p$. \label{rat} \end{thm}

\pf Let $\pi$ be the primary divisor of $p$ corresponding to the pair $(L,M)$ with both $L$~and~$M$ positive.  Then by Proposition~\ref{twothree}, Proposition~\ref{omega}, and Theorem~\ref{main}, the values of $\chi$ prescribed above agree with~$\chi_{\pi}$, so will consistently define a cubic character character modulo~$\pi$, and hence modulo~$p$.  Here we have used $\gamma=2+3\omega$ and $\gamma=4+3\omega$ in Theorem~\ref{main} to obtain rules (e) and (f), respectively.  The fact that $\chi(c)=1$ is equivalent to $c$ being a cubic residue has already been established. \fp

We remark that when $\chi(q)\ne 1$ then any of the last three rules may be employed; the theorem implies that any applicable rule will give the same value for~$\chi(q)$.  It is also possible that none of them can be used --- this occurs on average for one in every twenty-seven primes~$q$.  (Technically speaking, the above rules do not cover primes $q\equiv 2^{3m}5^{3n}\bmod 819$.  The smallest such prime is $q=181$.  The reader is invited to explain why the powers of $2^3$ and $5^3$ appear.)  In this case one would need to appeal to Theorem~\ref{main}, using a prime $l$ other than~7 or~13.

It is also interesting to note that because the rational functions in the last three rules give $\omega$ or $\bar{\omega}$-slopes, their numerators and denominators will never vanish at the primes $q\ge 5$ for which they apply.  (Due to Lehmer's result, the lines with slope 0 and $\infty$ contain points corresponding to primes for which $q$ is a cubic residue.)  However, as mentioned above, these functions duplicate the slopes in Lehmer's Theorem for the omitted values of~$q$, hence they vanish for three distinct values of~$t$, by Proposition~\ref{range}.  Applying this reasoning to rule (f), for example, provides an alternate method of verifying the following neat fact.

\begin{cor} The polynomial $t^3-7t-7$ factors completely modulo $q$ if $q\equiv\pm 1\bmod 7$, has a triple root when $q=7$, and otherwise is irreducible modulo~$q$. \end{cor}

\pf We have just argued that for primes $q\ge 5$ and $q\ne 7$, the denominator $t^3+27t^2-9t-27$ in rule (f) will vanish for three distinct values of $t$ modulo~$q$ when $q\equiv\pm 1\bmod 7$, but will never vanish otherwise.  In the latter case the polynomial must be irreducible, since any factorization would have to involve a linear factor.  This property is not affected by the substitution $t\mapsto 6t+9$, which yields the polynomial $6^3(t^3-7t-7)$.  By inspection the primes $q=2$, 3, and 7 agree with the statement, thus completing the argument. \fp

\noindent Thus the field extension $L/\mathbb{Q}$ defined by $t^3-7t-7$ (which turns out to be Galois) is the unique extension of~$\mathbb{Q}$ for which a prime $q\in\mathbb{Z}$ splits completely if and only if $q\equiv\pm1\bmod7$.  The prime $q=7$ ramifies, and all other primes are inert.

We conclude with a somewhat more exotic example: we shall determine whether $1982=2\cdot 991$ is a cubic residue modulo $p=\frac{1}{4}(3^{19}+5^{82})$, a 131-digit prime number.  One easily reads off $L=5^{41}$ and $M=3^8$.  Since $L\equiv M\equiv 1\bmod 4$, we find that $\chi(2)=\omega$ according to rule (a).  We also compute $991\equiv 1\bmod 9$, $991\equiv -3\bmod 7$, and $991\equiv 3\bmod 13$, so we may use either of rules (e) or (f), but not (d).  The relevant slope is
\[ \frac{M}{L} = \frac{3^8}{5^{41}} \equiv 672\equiv -319 \bmod 991. \]
We then use \textit{Pari-GP} to create a list of the values of the rational function in rule (e).  One discovers that 319 appears in the list (but not 672), so the negative sign must be taken, leading to $\chi(991)=\omega$.  For confirmation, we create a similar list of slopes using rule (f), and find that 672 appears this time rather than 319, so the positive sign holds.  Hence this rule also dictates that we should define $\chi(991)=\omega$.  In summary, \[ \chi(1982)=\chi(2)\chi(991)=\omega\cdot\omega\ne 1,\] so 1982 is in fact not a cubic residue modulo~$p$.

It is important to point out that the results presented here are primarily of theoretical rather than computational interest.  The \textit{Pari-GP} software invoked above can nearly instantly reach the same conclusion regarding the status of 1982 as a cubic residue modulo $p=\frac{1}{4}(3^{19}+5^{82})$ by simply reducing $1982^{(p-1)/3}$ modulo~$p$.

\section{Lehmer's Theorem}
Lehmer's original algorithm for determining the primes $p$ for which a given prime $q>3$ was a cubic residue was slightly more complicated than our version.  In summary, one writes $4p=L^2+27M^2$ with $L\equiv 1\bmod 3$, a restriction that can be removed.  If $q$ divides~$LM$, then $q$ is a cubic residue.  Otherwise, one first finds all quadratic residues $r$ modulo $q$ except for $r=1$.  One then finds the corresponding values of $u$ such that 
\[ r\equiv \frac{3u+1}{3u-3} \bmod q, \quad u\not\equiv 0, 1, -\frac{1}{2}, -\frac{1}{3}\bmod q. \]
Finally, one calculates $\mu^2$ satisfying
\[ \mu^2\equiv r\left(\frac{9}{2u+1}\right)^2 \bmod q. \]
Then the values of $\mu$ give the desired ratios $\frac{L}{M}$, that is,  $q$ is a cubic residue modulo $p$ if and only if $L^2\equiv \mu^2M^2\bmod q$.

To bring this into a nicer form, we first write $r\equiv t^2\not\equiv1\bmod q$, which gives the appropriate quadratic resides.  Solving for $u$ leads to
\[ u=\frac{3t^2+1}{3t^2-3} \quad\Longrightarrow\quad 2u+1 = \frac{9t^2-1}{3t^2-3}, \]
where the restrictions on $u$ require that $t\not\equiv 0, \pm\frac{1}{3},\pm\frac{1}{\sqrt{-3}}\bmod q$ in addition to $t\not\equiv\pm 1\bmod q$ from before.  Finally, we find that
\[ \mu \equiv \pm\frac{9t(3t^2-3)}{9t^2-1}\equiv\pm \frac{3t(3t+3)(3t-3)}{(3t+1)(3t-1)} \bmod q. \]
Since $\mu(-t)=-\mu(t)$ we may omit the $\pm$ sign.  Also, since $q>3$ we may replace $3t$ by $t$ and still obtain the same set of values for~$\mu$.  The final test for $q$ being a cubic residue then becomes
\[ LM\equiv0\bmod q \quad\mbox{or}\quad (t+1)(t-1)L\equiv t(t+3)(t-3) M\bmod q, \]
where $t\not\equiv 0,\pm 1, \pm 3, \pm\sqrt{-3}$.  However, the first several omitted values of~$t$ simply correspond to $L\equiv 0\bmod q$ or $M\equiv 0\bmod q$, so we combine the two statements into a single one and drop most of the restrictions on~$t$.  Finally, the values $t\equiv\pm\sqrt{-3}\bmod q$ correspond to slopes of $\frac{M}{L}\equiv\pm\frac{1}{3\sqrt{-3}}\bmod q$, the two slopes for which $L^2+27M^2$ is always divisible by~$q$.  So this qualification may be more clearly incorporated by simply requiring that $p\ne q$.  We have now reached the form given in equation~(\ref{lehmer}).

\end{document}